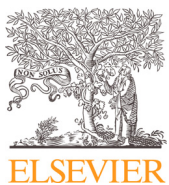
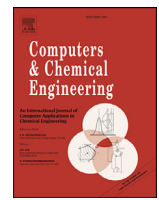

# Two-phase approaches to optimal model-based design of experiments: how many experiments and which ones?

Charlie Vanaret[a], Philipp Seufert[a], Jan Schwientek[a], Gleb Karpov[b], Gleb Ryzhakov[b], Ivan Oseledets[b,c], Norbert Asprion[d], Michael Bortz[a,*]

[a] *Fraunhofer Center for Machine Learning and Fraunhofer ITWM, Kaiserslautern, Germany*
[b] *Skoltech, Moscow, Russia*
[c] *INM RAS, Moscow, Russia*
[d] *BASF SE, Ludwigshafen am Rhein, Germany*

**A B S T R A C T**

Model-based experimental design is attracting increasing attention in chemical process engineering. Typically, an iterative procedure is pursued: an approximate model is devised, prescribed experiments are then performed and the resulting data is exploited to refine the model. To help to reduce the cost of trial-and-error approaches, strategies for model-based design of experiments suggest experimental points where the expected gain in information for the model is the largest. It requires the resolution of a large nonlinear, generally nonconvex, optimization problem, whose solution may greatly depend on the starting point. We present two discretization strategies that can assist the experimenter in setting the number of relevant experiments and performing an optimal selection, and we compare them against two pattern-based strategies that are independent of the problem. The validity of the approaches is demonstrated on an academic example and two test problems from chemical engineering including a vapor liquid equilibrium and reaction kinetics.

© 2021 The Authors. Published by Elsevier Ltd.
This is an open access article under the CC BY license (http://creativecommons.org/licenses/by/4.0/)

## 1. Motivation

Design of experiments (DoE) subsumes all methodologies for the systematic planning of experiments. The aim of DoE is to suggest a number of experiments as informative as possible, such that the parameters of a model may be estimated as reliably as possible. Ideally, forming the DoE optimization problem requires the knowledge of the model parameters, whose true values are unknown. In this paper we follow an iterative approach to estimate the model parameters $p$, starting from an initial guess. This is represented in the model validation and adjustment workflow in Fig. 1. Solving the DoE optimization problem provides an optimal design for the current $p$. The prescribed experiments are then carried out by the experimenter in order to gain information. With the observations, we update our estimate of $p$ and iterate until termination. The estimates of the model parameters are obtained by minimizing a loss function (the difference between the measured outputs of the system and the predictions of the model) ; it is usually a residual sum of squares, however other loss functions adapted to special situations, for example robust with respect to outliers, may also be used (Rousseeuw and Leroy, 2003).

Two questions arise in this context:

Q1. how many experiments are necessary to obtain accurate estimates of the model parameters?
Q2. how should the experiments be designed in order to maximize the reliability of the model based on estimates of the model parameters?

In the following, $P$ denotes the number of model parameters and $N$ the number of required experiments. Q1 concerns the minimal number of experiments $N$ required to estimate the $P$ parameters of the model. In the case of univariate observations, at least as many experiments as the number of model parameters must be performed: $N \geq P$. The model parameters may be estimated in the case $N = P$ if the experiments are chosen such that the covariance matrix of the estimates is not singular (Bates and Watts, 1988).

* Corresponding author.
*E-mail address:* michael.bortz@itwm.fraunhofer.de (M. Bortz).

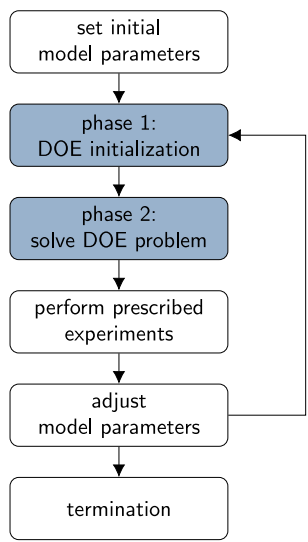

**Fig. 1.** Workflow for model validation and adjustment.

In the case of multivariate observations with $m$ outputs, estimating the model parameters requires at least $\lceil P/m \rceil$ experiments if the outputs are uncorrelated (the corresponding covariance matrix is not singular). In practice however, the $m$ outputs may stem from the model and may be significantly correlated. Estimating the model parameters may thus require at least $P$ experiments, depending on the properties of the covariance matrix (Fedorov and Leonov, 2014).

Q2 concerns the choice of the optimal experiments. Performing more than $P$ experiments increases the reliability of the estimates of the model parameters (Bates and Watts, 1988; Fedorov and Leonov, 2014) and consequently the predicted outputs of the model. Fedorov and Leonov (2014) showed that designing experiments such that the maximum error of the model-based predictions with respect to the experiments is minimal, is equivalent to minimizing the determinant of the covariance matrix of the model parameters. This optimization problem is generally nonlinear and nonconvex (Schöneberger et al., 2010). The quality of a local optimum depends on the initial point, the availability of exact derivatives and the choice of the optimization method (gradient descent, Newton-based methods, quasi-Newton methods, ...) Distinct initial points may lead to the same local optimum, and different optimization methods starting from the same initial point may reach distinct local optima. Furthermore, it is not possible to assess their global optimality. Zankin et al. (2018) made remarkable analytic progress for Vandermonde covariance matrices. Global optimization solvers certainly are an option, but for practical purposes, especially for fast proof-of-principle studies in industrial contexts, technical challenges may arise.

Fortunately, the approaches of statistical and linear experimental designs can be exploited for the initialization and multistarting of nonconvex experimental design optimization to increase the odds of finding the global optimum. Since Newton-based solvers converge quadratically in a neighborhood of a local optimum, finding accurate initial points is crucial for the overall convergence.

*Contributions* In this article, we describe four strategies with single start and multistart to initialize the nonconvex experimental design problem with satisfactory initial points and reduce the cost of the optimization process. We adopt a two-phase approach:

- phase 1: an initial point is generated:
  - by using a problem-independent pattern-based strategy (Section 3) in combination with single start and multistart, or
  - by solving a problem-dependent discretization of the experimental design problem that selects design points from a pool of candidates, approximates the optimal number of experiments and the values of the corresponding weights (Section 4).
- phase 2: the original experimental design problem is initialized with the phase-1 solution (Problem (3)).

When the number of experiments is optimal, the reliability of the model with respect to its parameters increases when the number of repetitions of some of the experiments is carefully tuned. However, in the nonlinear case, the questions of how many experiments to pick, which ones, and which experiments to repeat are non trivial and have not been addressed so far. This is one major contribution of this article. The algorithm MaxVol (Goreinov et al., 1997) originates from the family of algorithms derived for low-rank matrix approximations. It can quickly find a submatrix close to the one with the so-called D-optimality property (Mikhalev and Oseledets, 2018). The D-optimality criterion is widely used in the design of experiment practice (see the review of Hadigol and Doostan (2018)). To the best of our knowledge, MaxVol has not been used for experimental design in the literature so far. We propose a variant of MaxVol, wMaxVol, that accommodates multivariate outputs and approximates the weights of the relevant experiments. We also implemented a verification test to assess the optimality of the phase-1 and phase-2 solutions. It exploits the Kiefer-Wolfowitz equivalence theorem.

*Outline* We introduce the DoE problem and mathematical notations in Section 2. In Section 3, we describe two state-of-the-art pattern-based strategies, factorial design and quasi-random Sobol sequences, and discuss possible multistart strategies. We then introduce two discretization strategies in Section 4, the Weighted Discretization Approach inspired by Boyd and Vandenberghe (2004) and a variant of the MaxVol algorithm (Goreinov et al., 2010), that select the most relevant experiments among a fixed set of candidate experiments and determine their relevance. We introduce three test problems, an academic example and two chemical engineering use cases, in Section 5 and assess the benefits and limits of each strategy on the test problems in Section 6.

## 2. Model-based experimental design

Overviews of the DoE formalism and its application to chemical engineering problems can be found in Arellano-Garcia et al. (2007), Franceschini and Macchietto (2008). An exemplary case study on biodiesel production is given in Franceschini and Macchietto (2007). Mukkula and Paulen (2017) presented a generalization to DoE based on exact confidence intervals for highly nonlinear systems. The extension of the formalism to optimal control problems is discussed in Körkel (2002), the inclusion of stochastic uncertainties is the topic of Körkel et al. (2004). A joint consideration of both model discrimination and parameter estimates for time-dependent systems is given in Galvanin et al. (2016).

### 2.1. Notations

In model-based experimental design, the observations $y \in \mathbb{R}^m$ are given by:

$$y = f(x; p) + \varepsilon \qquad (1)$$

where:

- $x \in \mathcal{D} \subset \mathbb{R}^n$ is a vector of design variables of size $n$ ;
- $P$ is the number of model parameters ;
- $p \in \mathbb{R}^P$ is a vector of model parameters ;
- $f : \mathcal{D} \times \mathbb{R}^P \to \mathbb{R}^m$ is the model ;

- the observations are subject to a normally distributed error $\varepsilon \sim \mathcal{N}(0, \sigma^2)$ with zero mean and variance $\sigma^2$.

Following the definition from Fedorov and Leonov (2014), we call *experiment* the triplet $(X, r, y)$ where $X = (x_1, \ldots, x_N)$ is a collection of $N$ design points and $r_i$ is the number of repetitions of the design point $x_i$. The pair $(X, r)$ is called a *design of experiment*. We can also denote it as:

$$\xi = \begin{Bmatrix} X \\ w \end{Bmatrix} = \begin{Bmatrix} x_1 & \ldots & x_N \\ w_1 & \ldots & w_N \end{Bmatrix}, \quad (2)$$

where $w_i = \frac{r_i}{\sum_{i=1}^{N} r_i}$ are normalized weights. The points $(x_1, \ldots, x_N)$ are called the *spectrum* or the *support points* of the design and can be referred to as $supp(\xi)$.

From the viewpoint of a more general continuous design theory (Kiefer, 1959; Fedorov and Leonov, 2014), $w_i$ may vary continuously in [0, 1]. When the model produces a single output ($m = 1$) and is linear with respect to the model parameters $p$, each weight $w_i$ can be interpreted as the information gain or the importance of experiment $i$: experiments with low weights tend to be insignificant, while experiments with weights close to 1 tend to be meaningful. The magnitudes of the weights provide an indication of how sensitive the model is to a given design point under uncertainty: the number of measures of a new prescribed experiment $x_i$ should be proportional to the corresponding weight $w_i$ in order to reduce the uncertainty of the model. In continuous design, $\xi$ also denotes the probability measure on the domain $\mathcal{D}$ such that $\int_{\mathcal{D}} \xi(dx) = 1$; the reader should distinguish carefully between a design and its measure.

The *experimental design problem* consists in finding a design that minimizes a given statistical criterion $\Psi$:

$$\begin{aligned}
\min_{\xi := \{X, w\}^T} \quad & \Psi(\mathcal{I}(\xi)) \\
\text{s.t.} \quad & X := (x_1, \ldots, x_N) \in \mathcal{D}^N \\
& \sum_{i=1}^{N} w_i = 1 \\
& 0 \leq w_i \leq 1 \\
& 0 \leq c(x_i, y), \quad \forall i \in \{1, \ldots, N\}
\end{aligned} \quad (3)$$

where:

- $\mathcal{I}$ is the Fisher *information matrix*:

$$\mathcal{I}(\xi) = \int_{\mathcal{D}} \varphi(x) \varphi(x)^T \xi(dx) = \sum_{i=1}^{N} w_i \varphi(x_i) \varphi^T(x_i) = \sum_{i=1}^{N} w_i \mu(x_i) \quad (4)$$

where $\mu(x_i) = \varphi(x_i) \varphi^T(x_i)$ ;
- $\varphi : \mathbb{R}^n \to \mathbb{R}^{P \times m}$ is a matrix-valued function, following the standard notation in the literature (Fedorov, 2013). In nonlinear model-based experimental design, $\varphi(x_i)$ is taken as $J_p(x_i; t)^T$, the transpose of the Jacobian of the model $f$ at the point $x_i$ with respect to the model parameters;
- $0 \leq c(x_i, y), \forall i \in \{1, \ldots, N\}$ are constraints that may be imposed by the application that depend on $x$ and the observations $y$.

The traditional optimality criteria $\Psi$ are invariants of the information matrix $\mathcal{I}$: minimize the trace of $\mathcal{I}^{-1}$ (A-optimality), minimize the determinant of $\mathcal{I}^{-1}$ or equivalently maximize the determinant of $\mathcal{I}$ (D-optimality), maximize the minimum eigenvalue of $\mathcal{I}$ (E-optimality), maximize the trace of $\mathcal{I}$ (A*-optimality). In the rest of the paper, we focus exclusively on the log-D-criterion, that is the maximization of the decimal logarithm of the determinant of the information matrix (or equivalently, minimize the negative logarithm):

$$\Psi(\mathcal{I}(\xi)) = -\log \det(\mathcal{I}(\xi)) \quad (5)$$

This criterion results in maximizing the differential Shannon information content of the parameter estimates.

### 2.2. Parametric uncertainty

Under mild assumptions, the estimator of the model parameters converges to the true value in the model validation process used in this article (Fig. 1). A discussion and results can be found in (Fedorov and Leonov, 2014, p. 19). Alternatively, we can consider parametric uncertainty in the constraints and the information matrix. In this case, the uncertain model parameters $p$ are incorporated into the optimization problem. Traditional approaches include an average-case and a worst-case approach. The objective function and the constraints are replaced by the expectation with respect to $p$ or the maximum (minimum) with respect to $p$, respectively. We refer to (Fedorov and Leonov, 2014, pp. 80-84), Asprey and Macchietto (2002) and Körkel et al. (2004) for details.

### 2.3. Equivalence theorem

A significant result in the theory of optimal designs was the emergence of the so-called equivalence theorems for different optimality criteria, which established a connection between various formulations of optimization problems. A generalized modification of equivalence theorems for the multicriteria case can be found in (Fedorov and Leonov, 2014, pp. 68-69). Under mild assumptions on convexity and monotonicity of the statistical criterion $\Psi$ (Fedorov and Leonov (2014) showed that this is the case for the D-criterion), they show that for two designs $\xi$ and $\hat{\xi}$ with finite optimality criteria and for $\alpha \in (0, 1)$, there exists a function $\gamma$ such that:

$$\Psi\Big((1-\alpha)\mathcal{I}(\xi) + \alpha \mathcal{I}(\hat{\xi})\Big) = \Psi(\mathcal{I}(\xi)) + \alpha \int_{\mathcal{D}} \gamma(x, \xi) \hat{\xi}(dx) + e(\alpha; \xi, \hat{\xi}) \quad (6)$$

where $e$ is a function of $\alpha$ parameterized by $\xi, \hat{\xi}$, and $e(\alpha; \xi, \hat{\xi}) = o_{\alpha \to 0}(\alpha)$, that is $\lim_{\alpha \to 0} \frac{e(\alpha; \xi, \hat{\xi})}{\alpha} = 0$.

For the D-criterion, Fedorov and Leonov (2014) showed that:

$$\gamma(x, \xi) = P - d(x, \xi), \quad (7)$$

where $d$ is the *sensitivity function*:

$$d(x, \xi) = tr\big(\mathcal{I}(\xi)^{-1} \mu(x)\big) \quad (8)$$

The Kiefer-Wolfowitz equivalence theorem (Kiefer, 1959) states the following equivalent conditions:

1. the design $\xi^*$ minimizes $\Psi(\mathcal{I}(\xi))$ ;
2. the design $\xi^*$ minimizes $\max_{x \in \mathcal{D}} d(x, \xi)$ ;
3. the minimum over $\mathcal{D}$ of $\gamma(x, \xi^*)$ is 0, and occurs at the support points $x_i^* \in supp(\xi^*)$ of the design $\xi^*$.

Condition 3. can be written:

$$\begin{cases} d(x, \xi^*) \leq P, & \forall x \in \mathcal{D} \\ d(x_i^*, \xi^*) = P, & \forall x_i^* \in supp(\xi^*) \end{cases} \quad (9)$$

#### 2.3.1. Univariate observation ($m = 1$)

Using the invariance of the trace under cyclic permutations and the fact that $\varphi(x) \in \mathbb{R}^{1 \times P}$, we have:

$$d(x, \xi) = tr\big(\varphi^T(x) \mathcal{I}(\xi)^{-1} \varphi(x)\big) = \varphi^T(x) \mathcal{I}(\xi)^{-1} \varphi(x), \quad (10)$$

which coincides with the standardized variance of the predicted observation (Fedorov and Leonov, 2014; Atkinson, 2008). Eq. (9) thus simplifies to:

$$\begin{cases} \varphi^T(x) \mathcal{I}(\xi^*)^{-1} \varphi(x) \leq P, & \forall x \in \mathcal{D} \\ \varphi^T(x_i^*) \mathcal{I}(\xi^*)^{-1} \varphi(x_i^*) = P, & \forall x_i^* \in supp(\xi^*) \end{cases}$$

### 2.3.2. Multivariate observation ($m > 1$)

Often, the experimenter gathers $m > 1$ measurements simultaneously. We assume that $m$ components from one observation are correlated, but that distinct observations are independent. Unlike the univariate case, the information matrix $\mathcal{I}(\xi)$ includes the covariance matrix $\Sigma$:

$$\mathcal{I}(\xi) = \int_{\mathcal{D}} \varphi(x) \Sigma^{-1} \varphi^T(x) \xi(dx)$$
$$= \sum_{i=1}^{N} w_i \mu(x_i), \text{ where } \mu(x) = \varphi(x) \Sigma^{-1} \varphi^T(x). \quad (11)$$

Note that in the multivariate case, $1 \leq rank(\mu(x)) \leq m$, while in the univariate case (Eq. (4)), $rank(\mu(x)) = 1$.

## 3. Pattern-based strategies

### 3.1. Quasi-random Sobol sequences

#### 3.1.1. Construction

Sobol sequences (Sobol) are quasi-random low-discrepancy sequences of $p$ $n$-dimensional points (Sobol, 1967). They prevent the formation of clusters of samples points, and cover the domain more evenly than pseudorandom sequences (Kuipers and Niederreiter, 2012). The generation of Sobol sequences has been tremendously improved by Antonov and Saleev (1979), Jäckel (2002) over the years. A self-contained sequence of $p$ numbers in $(0, 1)$ is generated for each dimension, independently of the others. Sobol sequences have been successfully generated in very high dimension and have become a best practice for financial applications (Savine, 2018).

#### 3.1.2. Multistart approach

Common Sobol programming libraries implement Sobol sequences by generating sequences of integers between 0 and $2^{32} - 1$ in a recursive fashion. The first $n$ terms of a sequence in dimension $k > n$ exactly corresponds to the sequence in dimension $n$. An integer parameter `skip` of the Sobol generator controls the starting point in the sequence ; the first `skip` points are skipped. Calling the generator with different values of `skip` thus generates different Sobol sequences, which can be used as distinct initial points of the experimental design problem (Fig. 2). Skipped points are shown in light gray.

Note that consecutive integer values for `skip` produce Sobol sequences that differ by one point ; Fig. 2a shows that the gray point $(0.5, 0.5)$ is skipped and that the point $(0.96875, 0.59375)$ is added. More generally, two Sobol sequences of $p$ points generated with `skip` $= a$ and `skip` $= b > a$, respectively, will differ by $\min(b - a, p)$ points. For example, the Sobol sequence with `skip` $= 21$ (Fig. 2b) has no point in common with the Sobol sequence generated with `skip` $= 1$, since all 20 points were skipped.

### 3.2. Factorial designs

#### 3.2.1. Construction

A factorial design selects all (or a subset of) combinations of $n$ parameters at their lower and upper bounds. It studies the effects of individual parameters, as well as the effects of interactions between parameters. The **full** factorial design (Fig. 3a) generates all $S_f = 2^n$ combinations of lower and upper bounds in $n$ dimensions (the corners of the domain). A **reduced** or fractional factorial design (Fig. 3b) selects only a subset of combinations of the full factorial design and scales linearly with the dimension $n$ ; it generates $S_r = 2^{1+\lfloor \log_2(n) \rfloor} \approx 2n$ points. A reduced design avoids the redundancy of many experiments in a full design and describes the most important features of the problem.

Our strategy `Factorial` implements an extension of the factorial design to an arbitrary number of experiments $N$:

- if $N \leq S_r$, we select the first $N$ experiments of the reduced design ;
- if $S_r < N \leq S_f$, we select all $S_r$ experiments of the reduced design and the first $N - S_r$ of the remaining experiments of the full design ;
- if $S_f < N$, we select all $S_f$ experiments of the full design, and apply recursively `Factorial` on a subdomain of the domain (Fig. 3c) to select the remaining $N - S_f$. The subdomain may be chosen such that its volume is a given fraction (e.g. 50%) of the original domain.

#### 3.2.2. Multistart approach

The selection of experiments among the set of corners is inherently combinatorial, and lends itself well to multistart. For example, a design with $n = 3$ variables and $N = 6$ experiments may be generated by selecting the $S_r = 4$ experiments of the reduced design, then by taking the 2 remaining experiments among the $S_f - S_r = 8 - 4 = 4$ remaining corners of the full design. We thus have $\binom{4}{2} = 6$ possible combinations of experiments. To generate more experiments, we can freely select the 6 experiments among the 8 possible corners, that is $\binom{8}{6} = 28$ combinations. For a given number of multistart points, we shrink the domain recursively as many times as necessary, from which we can generate enough experiments.

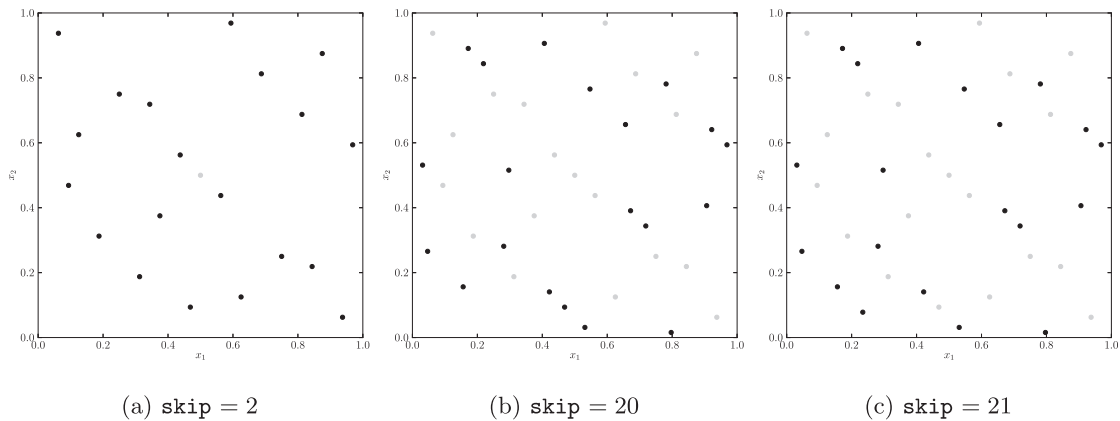

(a) `skip = 2`     (b) `skip = 20`     (c) `skip = 21`

**Fig. 2.** Sobol sequences with different `skip` values ($n = 2$, $p = 20$ points shown in black; the skipped points are shown in grey).

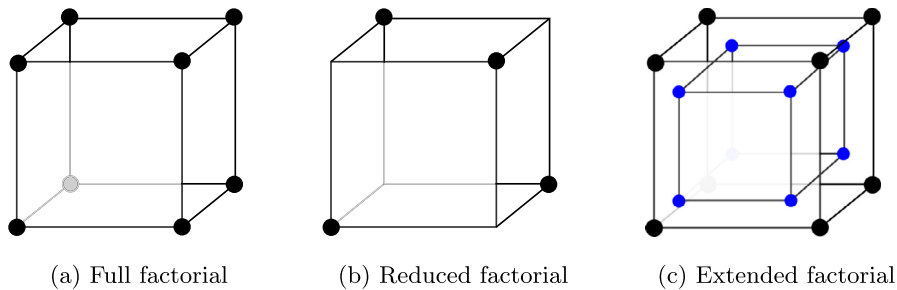

(a) Full factorial          (b) Reduced factorial          (c) Extended factorial

**Fig. 3.** Comparison of full, reduced and extended factorial designs ($n = 3$).

## 4. Discretization strategies

The experimental design problem on the continuous domain may be approximated by selecting experiments from a fixed finite set of $K$ candidate experiments. The resulting experiment selection problem exhibits mild properties (such as reduced size and convexity) and is thus tractable. It is an efficient option as phase 1 of the two-phase approach and is hoped to produce solutions close to good (possibly global) minimizers, thus reducing the computational effort of phase 2. This concept is similar to the two-phase approach in Dantzig's simplex method, in which the first phase aims at producing a feasible initial point.

In this section, we describe two discretization strategies that solve the experiment selection problem:

- the Weighted Discretization Approach (WDA) solves a continuous optimization problem with respect to the $K$ weights and provides the number of relevant experiments (the number of positive weights) and their relevance (the magnitudes of the weights) ;
- the binary formulation (wMaxVol) is a greedy algorithm that determines a satisfactory combination of candidate experiments and the corresponding weights.

The experiment selection problem exhibits favorable practical assets. Since only the weights $w$ are adjusted, the Jacobian and Fisher information matrices can be computed once and for all for each candidate experiment. The discretization $S$ can be refined adaptively, should the simulations be costly.

The number of optimal experiments is approximated by the discretization-based strategies as a by-product of the experiment selection problem: it corresponds to the number of nonzero weights at the optimal solution. Should the number of experiments computed in phase 1 be too large an upper bound of the optimal number of experiments, certain weights may be set to 0 by the optimization solver to disable the corresponding design points.

### 4.1. Parameterizations

Since WDA and wMaxVol can be parameterized by a finite set $S$ of candidate experiments, we adopt the functional notations WDA(S) and wMaxVol(S). $S$ can be instantiated in an arbitrary manner, for instance:

- Grid: experiments are positioned according to a fixed rectangular topology ;
- Sobol and Factorial (Section 3) ;
- the user can include preferred experiments to increase the confidence of their expertise, and select particular subsets of experiments to rerun the optimization process ;
- a phase-2-optimal solution achieved with another strategy or itself (in a recursive fashion).

and can be composed with a filtering strategy that discards experiments that do not satisfy a set of constraints (Feasible(S)).

The size of the candidate set may also be specified, for example Grid(50) or Sobol(2).

A possible globalization approach for a discretization strategy is to perform multistart with several parameterizations, e.g. WDA(Grid), WDA(Feasible(Grid)) and WDA(Sobol).

### 4.2. Weighted discretization approach (WDA)

Boyd and Vandenberghe (2004) suggested to approximate the experimental design problem (Eq. (3)) on a fixed set $\bar{X}$ of $K$ candidate experiments. The size of the problem is thus reduced from $N \times (n+1)$ variables (the experiments and their weights) to $K$ variables (the weights). The number of positive weights of the optimal solution determines the number of relevant experiments, while the magnitude of the weights indicate the relevance of the optimal experiments. The Weighted Discretization Approach (WDA) problem is:

$$\begin{aligned}\min_{w \in \mathbb{R}^K} \quad & \Psi\left(\mathcal{I}\left(\left\{\begin{array}{c}\bar{X}\\w\end{array}\right\}\right)\right)\\ \text{s.t.} \quad & \sum_{i=1}^K w_i = 1\\ & 0 \leq w_i \leq 1\end{aligned} \quad (12)$$

Its minimum is an upper bound of the minimum of the experimental design problem.

The WDA problem is a positive semi-definite optimization problem with respect to the weights $w$ that is convex (minimization of a convex function subject to convex constraints) when $\Psi$ is the A-criterion, the D-criterion or the E-criterion. In this favorable case, the optimization problem can be solved to global optimality using state-of-the-art convex optimization techniques that can generally handle large-scale problems. Since the analytical derivatives of the objective function and the constraints with respect to the weights $w$ are known analytically, robust off-the-shelf NLP solvers (such as IPOPT, SNOPT or MINOS) are also suited.

### 4.3. wMaxVol algorithm

MaxVol and its successor rect_MaxVol are greedy algorithms that originate from the family of low-rank matrix approximation algorithms (Goreinov and Tyrtyshnikov, 2001; Goreinov et al., 1997). The most general approach, rect_MaxVol, is looking for a quasi-optimal solution to the following optimization problem. Consider a tall matrix $A \in \mathbb{R}^{n \times m}$. The problem is to find the submatrix $A(I)$ composed of rows of $A$ enumerated in $I \subset \mathcal{P} = \{1, \ldots, n\}$, such that $|I| = k$ and $m \leq k < n$, solution to:

$$\max_{I, |I|=k} \begin{cases} |\det(A(I))| & \text{if } k = m \\ \det\left(A(I)^T A(I)\right) & \text{if } k > m \end{cases} \quad (13)$$

#### 4.3.1. Rectangular MaxVol

The general algorithm rect_MaxVol (Algorithm 1) consists in two steps:

**Algorithm 1:** rect_MaxVol.

**Data**: full-rank matrix $A \in \mathbb{R}^{n \times m}$, $k \in \mathbb{N}$, $n > k \geq m$, tolerance $\delta \geq 0$
**Result**: $k \times m$ dominant submatrix $\tilde{A}$
1  $(\tilde{A}, \tilde{\mathcal{P}}) \leftarrow \text{MaxVol}(A, \delta)$ (Goreinov et~al., 2010) ;
2  $\mathcal{P} \leftarrow \tilde{\mathcal{P}}(1:m)$ ;
3  **while** $|\mathcal{P}| \leq k$ **do**
4  $\quad C \leftarrow A\tilde{A}^+$ ;
5  $\quad i \leftarrow \arg\max_{i \in \tilde{\mathcal{P}} \setminus \mathcal{P}} \|C(i,:)\|_2^2$ ;
6  $\quad \mathcal{P} \leftarrow \mathcal{P} \cup \{\tilde{\mathcal{P}}(i)\}$ ;
7  $\quad \tilde{A} \leftarrow \begin{bmatrix} \tilde{A} \\ A(i,:) \end{bmatrix}$ ;
8  **end**
9  **return** $\tilde{A}$ ;

**Algorithm 2:** wMaxVol.

**Data**: $K = |\mathcal{D}|$, full-rank design matrix $A \in \mathbb{R}^{Km \times P}$, $n_{iter}$, $m \in \mathbb{N}^+$
**Result**: indices of distinct design points $\mathcal{P}_U$, corresponding weights $w$
1  $l \leftarrow \lfloor \frac{P}{m} \rfloor + 1$ ;
2  Randomly generate $\mathcal{L} \subset \{1, \ldots, K\}$ such that $|\mathcal{L}| = l$ and $A(\mathcal{L})^T A(\mathcal{L})$ is not singular ;
3  $\tilde{A} \leftarrow A(\mathcal{L})$ ;
4  $\mathcal{P} \leftarrow \{\}$ ;
5  **for** $k \in \{0, \ldots, n_{iter}\}$ **do**
6  $\quad C \leftarrow A\tilde{A}^+$ ;
7  $\quad i \leftarrow \arg\max_{i \in \{1, \ldots, K\}} \det\left(I + C_i C_i^T\right)$ ;
8  $\quad \mathcal{P} \leftarrow \mathcal{P} \cup \{i\}$ ;
9  $\quad \tilde{A} \leftarrow \begin{bmatrix} \tilde{A} \\ A_i \end{bmatrix}$ ;
10 **end**
11 $(\mathcal{P}_U, \mathcal{W}) \leftarrow$ (unique elements of $\mathcal{P}$, number of repetitions) ;
12 $\mathcal{W}(\mathcal{L} \cap \mathcal{P}_U) \leftarrow \mathcal{W}(\mathcal{L} \cap \mathcal{P}_U) + 1$ ;
13 $w \leftarrow \frac{\mathcal{W}}{n_{iter}}$ (weight normalization) ;
14 **return** $\mathcal{P}_U, w$ ;

1. the algorithm MaxVol (line 1): extraction of a square submatrix $\tilde{A}$ from $A$ by iteratively letting one row of $A$ enter the submatrix and another exit it (Goreinov et al., 2010) ; the greedy row swap is performed such that the determinant of the submatrix is maximized. $\tilde{A}$ is a quasi-optimal solution to Eq. (13). The algorithm terminates when no additional row swap can increase the objective function by more than $1 + \delta$ times ($\delta \geq 0$) or when the iteration limit is reached. The algorithm returns $\tilde{A}$ and $\tilde{\mathcal{P}}$, a vector of row permutations of $A$ such that the submatrix $\tilde{A}$ is composed of the first rows indexed by $\tilde{\mathcal{P}}$.
2. the greedy expansion of the square submatrix until it reaches the required shape (Mikhalev and Oseledets, 2018).

Matrix inverses, pseudo-inverses and matrix-matrix multiplications (lines 1 and 4 in Algorithm 1) are bottlenecks that may dramatically slow down the computations for relatively large matrices. Goreinov et al. (2010), Mikhalev and Oseledets (2018) introduced computationally efficient update formulas based on linear algebra for obtaining intermediate matrices at each iteration ; in practice, no costly inversion or multiplication is thus performed.

This problem is similar to the the D-optimal experimental design problem (Eqs. (3) and (5)) for an input matrix $A$ constructed in a way such that its rows are vector functions $\varphi(x)^T = J_p(x; p) \in \mathbb{R}^{1 \times P}$ computed for all $x \in \mathcal{D} \subset \mathbb{R}^n$. The submatrix sought by MaxVol is the design matrix whose rows are then the respective design points, which demonstrates the relevance of MaxVol for solving experimental design problems. However, the experimental design problem does not fit effortlessly within the off-the-shelf MaxVol framework for several reasons:

- in the standard framework, the number of rows $k$ in the resulting submatrix is a hyperparameter. In experimental design, $k$ corresponds to the number of points in the design and is an unknown quantity ;
- MaxVol does not produce weights. Applying MaxVol to experimental design thus results in equal-weighted designs ;
- in multivariate experimental design, the block $\varphi(x)^T \in \mathbb{R}^{m \times P}$ corresponding to the design point $x$ is composed of $m$ dependent rows. By analogy with the univariate case (the design matrix $A$ consists of $N$ rows), $A$ should consist of $N$ blocks. However, this block structure is not supported by the MaxVol framework as it operates on each row of $A$ separately. Therefore, it cannot produce a consistent design.

These reasons motivated the development of a MaxVol-based algorithm that corrects the weaknesses of the original algorithm and supports the resolution of experimental design.

*4.3.2. A variant of MaxVol for D-optimal experimental design*

We propose a variant wMaxVol (Algorithm 2) of rect_MaxVol for experimental design. To handle multivariate outputs, we introduce a block matrix form

$$A = \begin{pmatrix} A_1 \\ \vdots \\ A_K \end{pmatrix}, \tag{14}$$

where each block $A_i \in \mathbb{R}^{m \times P}$ denotes a separate single entity. In the case of experimental design, $A$ corresponds to the vertical concatenation of $K$ blocks $A_i := \phi(x_i)^T$ for each design point $x_i \in \mathcal{D} \subset \mathcal{R}^n$.

The row swap operation and the single-row expansion were replaced by their block counterparts, and the corresponding update formulas were derived to preserve the numerical efficiency of the algorithm. Here, $i$ denotes a block index. wMaxVol can approximate the weights of the experiments by performing additional iterations ; the repetition of choices was implemented by picking $i \in \tilde{\mathcal{P}}$ instead of $i \in \tilde{\mathcal{P}} \setminus \mathcal{P}$ on line 5 of Algorithm 1. The constraint on the required size of the resulting submatrix was also dropped. The weights of the distinct experiments are obtained by normalizing the number of occurrences of the experiments. wMaxVol is therefore suited for approximately solving D-optimal experimental design, which we formally prove in Appendix A for the univariate case. The proof for the multivariate case is beyond the scope of this paper and will be given in a separate article. We have however performed a statistical study in Appendix B to analyze the convergence of wMaxVol on a multivariate example (Chebyshev polynomials and their derivatives). The analysis is supported by the Kiefer-Wolfowitz equivalence theorem (Section 2.3).

## 5. Test problems

In this section, we describe three test problems with increasing complexity on which the validity of our strategies will be assessed: an academic exponential example whose solution is known analytically, a flash distillation example governed by the MESH equations and an inequality-constrained tubular reactor example.

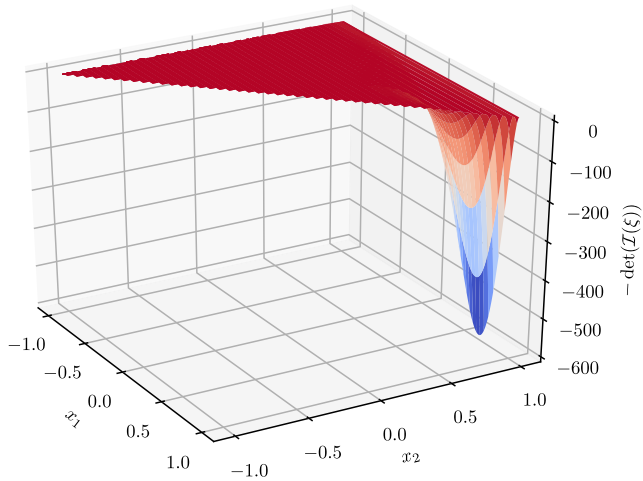

**Fig. 4.** Representation of the exponential problem (Eq. (19)) for $p = (1, 3)$ and $w = (0.5, 0.5)$.

## 5.1. Exponential example

We introduce a model represented by a one-dimensional function of input $x \in [-1, 1]$ and parameterized by $p = (p_1, p_2) > 0$:

$$f(x; p) = p_1 e^{p_2 x} \quad (15)$$

whose Jacobian with respect to $p$ is:

$$J_p(x; p) = \begin{pmatrix} e^{p_2 x} & p_1 x e^{p_2 x} \end{pmatrix} \quad (16)$$

Let $\xi = \begin{Bmatrix} x_1 & x_2 \\ w_1 & w_2 \end{Bmatrix}$ be a design composed of two experiments $x_i \in [-1, 1]$ and their corresponding weights $w_i \in [0, 1]$. The Fisher matrix $\mathcal{I}(\xi)$ can be written:

$$\mathcal{I}(\xi) = \sum_{i=1}^{2} w_i J_p(x_i; p)^T J_p(x_i; p)$$
$$= \begin{pmatrix} w_1 e^{2p_2 x_1} + w_2 e^{2p_2 x_2} & p_1(w_1 x_1 e^{2p_2 x_1} + w_2 x_2 e^{2p_2 x_2}) \\ p_1(w_1 x_1 e^{2p_2 x_1} + w_2 x_2 e^{2p_2 x_2}) & p_1^2(w_1 x_1^2 e^{2p_2 x_1} + w_2 x_2^2 e^{2p_2 x_2}) \end{pmatrix} \quad (17)$$

Its determinant is given by:

$$\det(\mathcal{I}(\xi)) = w_1 w_2 p_1^2 (x_1 - x_2)^2 e^{2p_2(x_1 + x_2)} \quad (18)$$

which is symmetric in $x_1$ and $x_2$. The D-optimal design $\xi^*$ is therefore the solution to the following optimization problem:

$$\begin{aligned} \min_{X \in \mathbb{R}^2, w \in \mathbb{R}^2} & \quad -w_1 w_2 p_1^2 (x_1 - x_2)^2 e^{2p_2(x_1 + x_2)} \\ \text{s.t.} & \quad -1 \leq X \leq 1 \\ & \quad 0 \leq w \leq 1 \\ & \quad x_1 \leq x_2 \\ & \quad w_1 + w_2 = 1 \end{aligned} \quad (19)$$

where the constraint $x_1 \leq x_2$ is added to break symmetries. Since the line $x_1 = x_2$ is a continuum of local maximizers for which $\det(\mathcal{I}(\xi)) = 0$, this constraint is not active at the solution (a minimum) of the problem.

In the rest of the paper, we set the parameters to $p = (1, 3)$ (Fig. 4). The method of Lagrange multipliers yields the global minimizer $\xi^* = \begin{Bmatrix} X^* \\ w^* \end{Bmatrix} = \begin{Bmatrix} \frac{2}{3} & 1 \\ 0.5 & 0.5 \end{Bmatrix}$ and the global minimum $-\det(\mathcal{I}(\xi^*)) = -\frac{e^{10}}{36}$.

The sensitivity function $x \mapsto d(x, \xi^*)$ is given by:

$$d(x, \xi^*) = \varphi^T(x) \mathcal{I}(\xi^*)^{-1} \varphi(x)$$

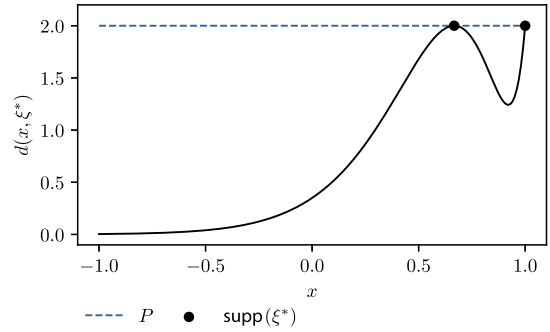

**Fig. 5.** D-optimality for the exponential example with $p = (1, 3)$.

$$= \left(18 x^2 (e^2 + 1) - 12 x (3 e^2 + 2) + 2(9 e^2 + 4)\right) e^{6x - 6} \quad (20)$$

and is represented in Fig. 5. The optimality bound $P = 2$ is reached at the support points of the D-optimal design $\xi^*$.

### 5.2. Chemical engineering problem: Flash distillation

The flash distillation example (Asprion et al., 2019) is represented in Fig. 6. A liquid mixture, called the feed, composed of methanol and water enters the unit with fixed flow rate $F$ and compositions $(z_m, z_w)$ with $z_m + z_w = 1$. The mixture is heated up by an external heat source with heat duty $\dot{Q}$ and partially vaporized. The produced liquid and vapor are at equilibrium at pressure $P$ and temperature $T$. Their respective liquid $(x_m, x_w)$ and vapor $(y_m, y_w)$ compositions depend on the degree of vaporization.

The flash unit enforces the so-called MESH equations:

- Mass balances

$$F z_m = V y_m + L x_m \quad (21)$$

$$F z_w = V y_w + L x_w \quad (22)$$

- Equilibrium

$$P y_m = P_m^0(T) x_m \gamma_m(x, T) \quad (23)$$

$$P y_w = P_w^0(T) x_w \gamma_w(x, T) \quad (24)$$

- Summation

$$x_m + x_w = y_m + y_w = z_m + z_w = 1 \quad (25)$$

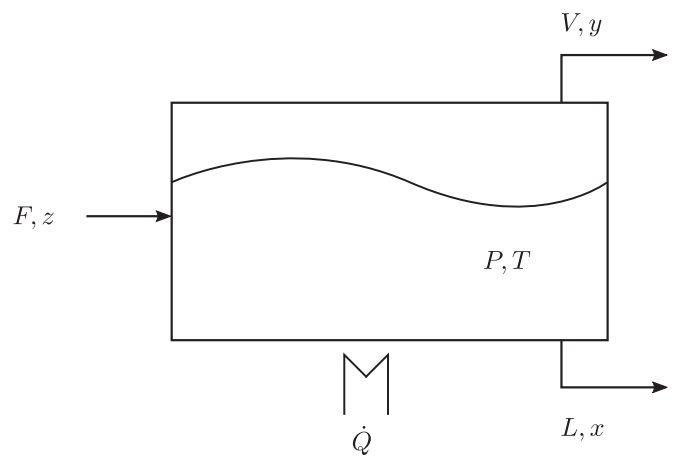

**Fig. 6.** Flash distillation unit.

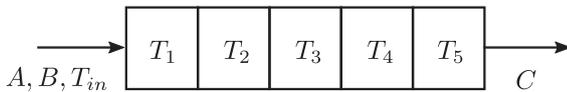

Fig. 7. Tubular reactor.

- Heat balance

$$\dot{Q} + FH_L(z, T_F) = VH_V(y, T) + LH_L(x, T) \quad (26)$$

where the vapor pressure of the pure elements $P^0$, the activity coefficients of the components in the mixture $\gamma$ and the enthalpies of the liquid and the vapor streams $H_L$ and $H_V$ are given by thermodynamic models. Additionally, the total molar flow $F$ of the feed is set to 1 kmol/h and the vapor molar flow $V$ is set to $10^{-6}$ kmol/h, which ensures a measurement at the boiling point of the mixture. This results in 10 equations and 12 unknowns ($L$, $V$, $F$, $x_m$, $x_w$, $y_m$, $y_w$, $z_m$, $z_w$, $P$, $T$, $\dot{Q}$). The experimental design problem has the following representation:

- 2 design variables: the pressure $P \in [0.5, 5]$ bar at equilibrium and the molar concentration of methanol $z_m \in [0, 1]$ of the feed. This amounts to fixing the last two degrees of freedom of the model ;
- 4 model parameters: the parameters $A_{12} = 0.3, A_{21} = 5.166, B_{12} = 1118$ and $B_{21} = -1473$ for the activity coefficients $\gamma_m$ and $\gamma_w$ of the non-random two-liquid (NRTL) model $\gamma$ ;
- 2 outputs: the concentration of methanol $y_m$ of the vapor output stream and the temperature $T$ at equilibrium.

The inverse of the covariance matrix is a diagonal matrix $\Sigma^{-1} = diag(10^{-2}, 10^4)$. The globally optimal solution to the experimental design problem for the flash distillation is unknown.

### 5.3. Chemical engineering problem: Tubular reactor

The tubular reactor example is an inequality-constrained problem that models the reaction:

$$A + 2B \rightarrow 3C \quad (27)$$

where A and B represent two reactants and C one product. This reaction takes place in a tubular reactor consisting of 5 sections (Fig. 7). A feed containing reactants A and B enters the reactor. After a partial conversion, a stream consisting of the product C as well as unconverted reactants exits the reactor.

The experimental design problem has the following representation:

- 2 design variables: $F_{ratio} = B/A \in [0.2, 1]$ is the feed ratio between the two molar flow rates $A$ and $B$, and $T_{in} \in [70, 110]°C$ is the temperature of the feed ;
- 2 model parameters: $E_A = 23500$ kJ/kmol is the activation energy of the reaction and $k_0 = 10^4$ is a pre-exponential factor in the expression of the kinetic rate $r = k x_A x_B^2$, where $k = k_0 e^{-\frac{E_A}{RT}}$ is the rate constant ;
- 8 outputs: $x_A$, $x_B$ and $x_C$ are the concentrations of the chemical substances in the output stream and $\{T_1, \ldots, T_5\}$ are the temperatures in the sections of the tubular reactor ;
- 3 inequality constraints: the concentrations $x_A$ and $x_B$ should be positive and the maximal temperature in the tubular reactor should be below $125°$.

The inverse of the covariance matrix is a diagonal matrix $\Sigma^{-1} = diag(10^4, 10^4, 10^4, 10^{-2}, 10^{-2}, 10^{-2}, 10^{-2}, 10^{-2})$. The globally optimal solution to the experimental design problem for the tubular reactor is unknown.

## 6. Numerical results

In this section, we assess the validity of the `Sobol`, `Factorial`, `WDA` and `wMaxVol` strategies on the test problems presented in Section 5. We detail here the numerical methods that were used:

- the simulations are solved with a Newton-based method within CHEMASIM, the BASF in-house program (Asprion et al., 2015), on an Intel i5 CPU @ 2.67GHz with 4GB RAM ;
- the Jacobian matrices are computed using analytical derivatives (CHEMASIM) ;
- the phase-1 WDA problem is solved using the SCS solver from the Python library `cvxpy` with a tolerance of $10^{-3}$ ;
- the phase-1 wMaxVol problem is solved using a Python implementation with a tolerance of $10^{-3}$. 1,000 iterations were performed to compute an estimate of the weights ;
- the phase-2 problem is solved using the NLP solver NLPQLP (Schittkowski, 2006) with a tolerance of $10^{-6}$ (CHEMASIM) ;
- a small multiple of the identity matrix (with a factor $10^{-8}$) was added as a regularization term to the phase-2 information matrix to avoid singularity (see Asprion et al. (2019));
- the matrix $\varphi^T(x) = J_p(x; p)$ is scaled with the diagonal matrix $diag(p) \in \mathbb{R}^{P \times P}$, where $p$ are the current model parameters. This results in a scaling of the information matrix $\mathcal{I}$.

### 6.1. Exponential example

WDA and wMaxVol produce the same solution on an unidimensional `Grid` design with $N = 11$ equidistant candidate experiments in $[-1, 1]$ (Fig. 8a). The model parameters are $p = (1, 3)$, for which the D-optimal design is $\xi^* = \begin{Bmatrix} \frac{2}{3} & 1 \\ 0.5 & 0.5 \end{Bmatrix}$. The nonzero weights are represented by dashed (phase 1) or colored (phase 2) disks centered in the corresponding candidates, and whose areas are proportional to the weights. The phase-1 optimal design is $\xi_{11} = \begin{Bmatrix} 0.6 & 1 \\ 0.5 & 0.5 \end{Bmatrix}$, close to the optimum. The phase-2 solution corresponds to the D-optimal design. Fig. 8c illustrates the sensitivity function $x \mapsto d(x, \xi_{11})$, along with the optimality bound $P = 2$ and the evaluation of all 11 candidates. It validates the global optimality of $\xi_{11}$ over the discrete grid, since the first experiment $x = \frac{2}{3}$ of $\xi^*$ cannot be achieved on the grid points.

The grid can be subsequently refined at will ; in Fig. 8b, an additional candidate experiment is added at $x = 0.7333$. The phase-1 optimal design becomes $\xi_{12} = \begin{Bmatrix} 0.6 & 0.7333 & 1 \\ 0.37 & 0.13 & 0.5 \end{Bmatrix}$: the experiment $x = 0.6$ with weight 0.5 is now split between $x = 0.6$ (weight 0.37) and $x = 0.7333$ (weight 0.13). During phase 2, the weight of the experiment $x = 0.7333$ is set to 0 by the solver, thus disabling it. The experiment $x = 0.6$ is refined to produce an experiment $x = \frac{2}{3}$ with weight 0.5. The optimality verification is shown in Fig. 8d: it certifies that $\xi_{12}$ is optimal on the discrete domain.

### 6.2. Flash distillation

#### 6.2.1. Single start

Fig. 9a and b portray the phase-1 (dashed) and phase-2 (colored) solutions for `Sobol` and `Factorial` with $N = 6$ experiments. The initial phase-1 weights are $\frac{1}{N}$. The corresponding phase-2 solutions are 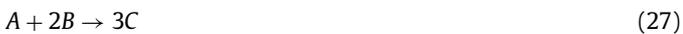 and

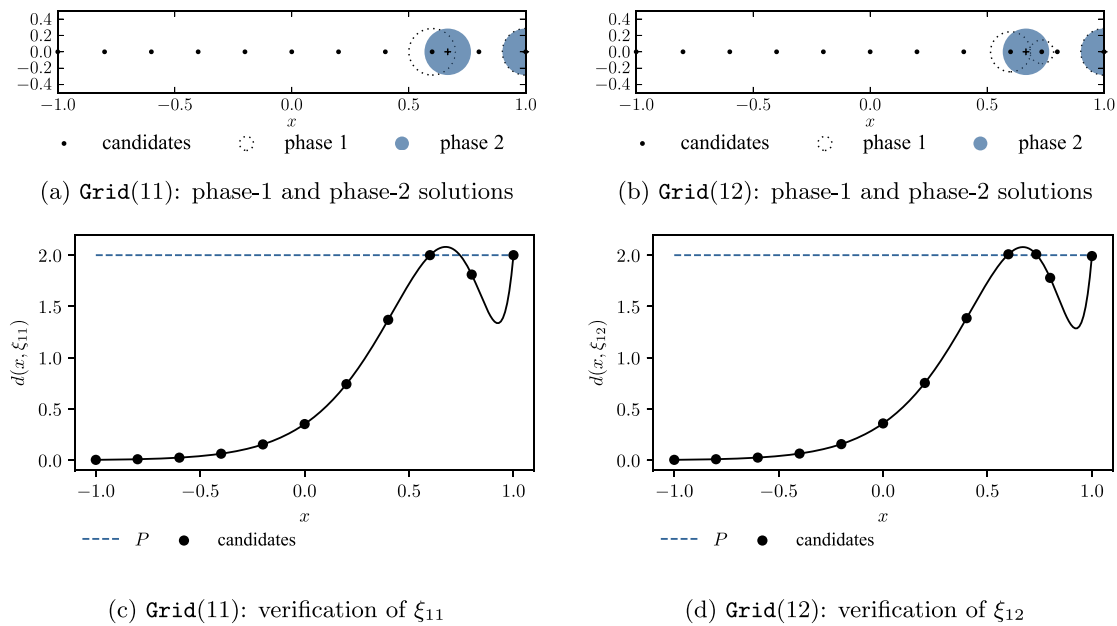

Fig. 8. Exponential example.

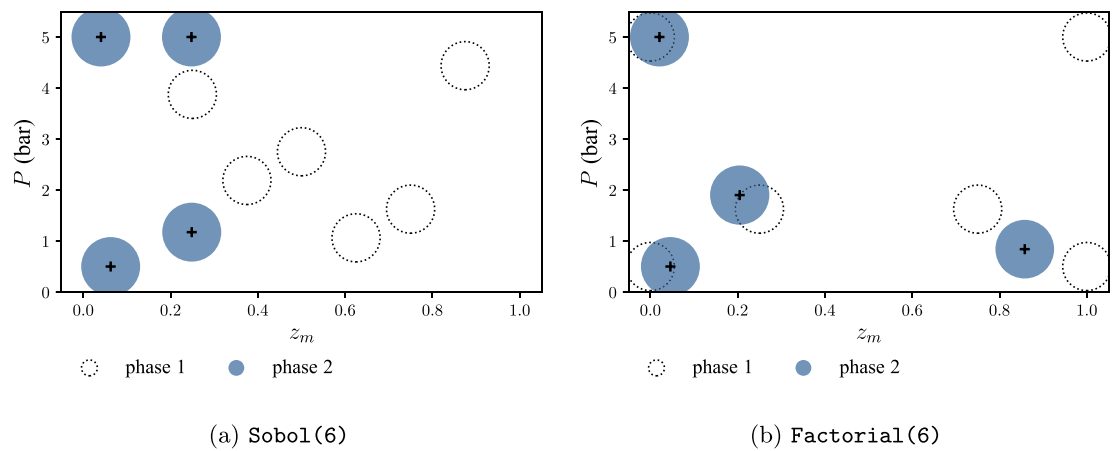

Fig. 9. Sobol and Factorial designs for the flash distillation with 6 experiments.

$$\begin{Bmatrix} 0.046 & 0.021 & 0.205 & 0.858 \\ 0.5 & 5 & 1.900 & 0.841 \\ 0.250 & 0.250 & 0.252 & 0.248 \end{Bmatrix},$$ respectively. Two of the experiments vanish in phase 2 as their weights are driven to 0 by the optimizer.

Fig. 10a illustrates the phase-1 and phase-2 solutions of the discretization strategies on a coarse grid of 25 candidate experiments. Extremes values (0 and 1) of $z_m$ correspond to pure components. In this case, $P$ and $T$ are directly coupled and it is not possible to determine the interaction parameters from the sole component data. Consequently, the grid is generated with $z_m$ ranging from 0.1 to 0.9. WDA and wMaxVol produce the same phase-1 solution with 6 nonzero weights, and therefore the same phase-2 solution $\hat{\xi}$. The optimality verification is shown in Fig. 10b. The candidates are located at the middle of each tile, and the color amplitude varies with the magnitude of the sensitivity $d$. Since the maximal value of $d$ is $P = 4$ at the phase-1 solution, the global optimality of the phase-1 solution is certified on the discrete domain.

The corresponding optimal phase-1 and phase-2 log-D-criteria $\Phi$ and CPU times are gathered in Table 1. For the discretization strategies, the CPU times include i) the generation of Jacobians at the candidate experiments (0.4s) and ii) the convergence time of phase 1. It shows that Factorial(6) achieves a local minimum whose evaluation is worse than the local minimum achieved by Sobol(6), WDA(Grid(25)) and wMaxVol(Grid(25)). This corroborates the fact that the experimental design problem is multimodal. wMaxVol proves significantly faster than WDA ; using a different convex solver may reduce the gap between WDA and wMaxVol.

Table 1
Optimal log-D-criteria and CPU times for the flash distillation ($N = 6$) with single start.

| Strategy | | $\Phi$ | CPU time (s) |
| --- | --- | --- | --- |
| Sobol(6) | Phase 1 | −3.263 | 25 |
| | Phase 2 | **−7.928** | |
| Factorial(6) | Phase 1 | 4.365 | 26 |
| | Phase 2 | −5.325 | |
| WDA(Grid(25)) | Phase 1 | −6.472 | 0.4 + 191 |
| | Phase 2 | **−7.928** | 15.13 |
| wMaxVol(Grid(25)) | Phase 1 | −6.472 | 0.4 + 0.312 |
| | Phase 2 | **−7.928** | 15.13 |

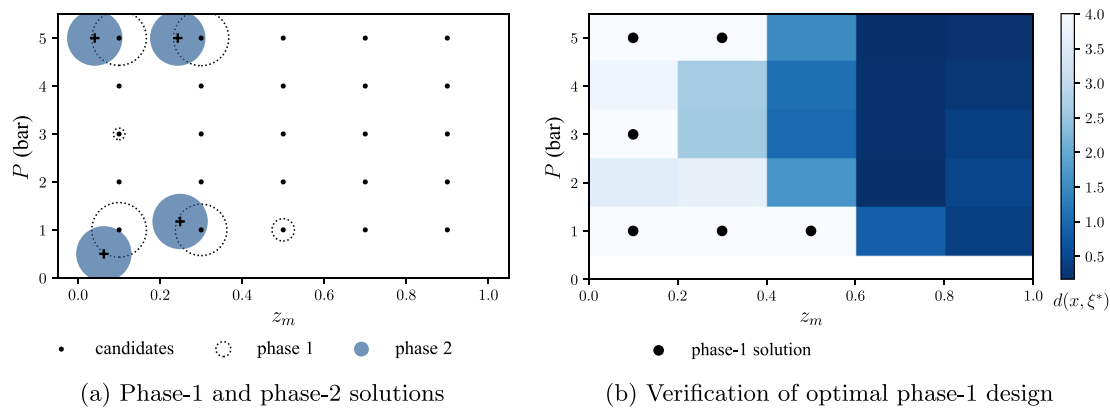

**Fig. 10.** WDA and wMaxVol designs for the flash distillation with 25 candidates.
(a) Phase-1 and phase-2 solutions   (b) Verification of optimal phase-1 design

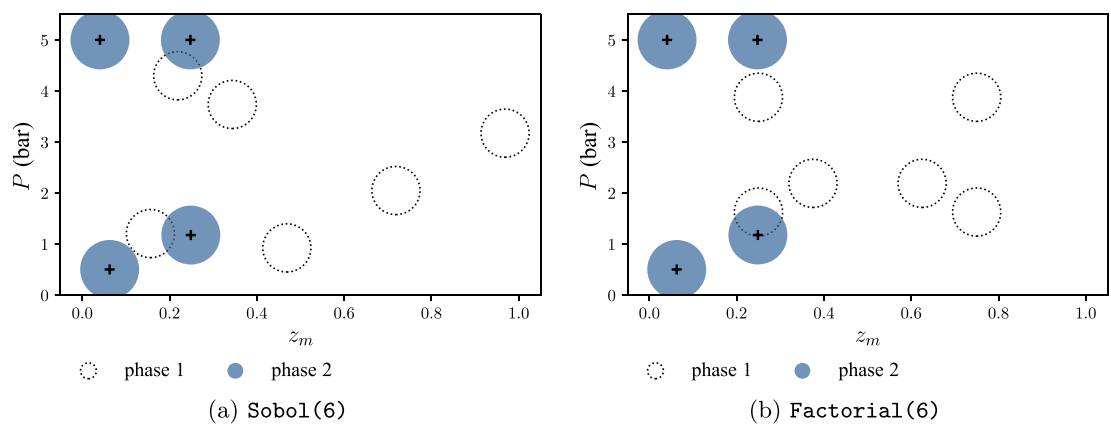

**Fig. 11.** Sobol and Factorial designs for the flash distillation with 6 experiments and 5 multistart runs.
(a) Sobol(6)   (b) Factorial(6)

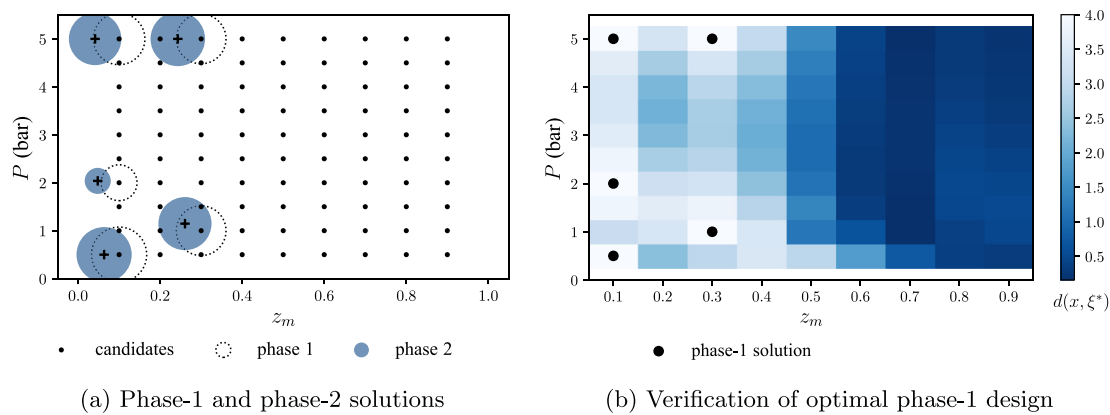

**Fig. 12.** WDA and wMaxVol designs for the flash distillation with 90 candidates.
(a) Phase-1 and phase-2 solutions   (b) Verification of optimal phase-1 design

*6.2.2. Multistart*

Fig. 11a and b illustrate the phase-1 and phase-2 solutions of Sobol and Factorial for the best run among 5 multistart runs. Both strategies produce the local minimum $\hat{\xi}$ (two weights are driven to 0 during phase 2).

Fig. 12a illustrates the phase-1 and phase-2 solutions of the discretization strategies on a fine grid ($9 \times 10$ points). WDA and wMaxVol produce the same phase-1 solution with 5 nonzero weights, very close to the putative global minimizer

$$\xi^* = \begin{Bmatrix} 0.048 & 0.042 & 0.063 & 0.261 & 0.243 \\ 2.039 & 5 & 0.5 & 1.147 & 5 \\ 0.055 & 0.224 & 0.248 & 0.230 & 0.242 \end{Bmatrix},$$

which is subsequently reached during phase 2. The global optimality of the phase-1 solution is certified on the discrete domain (Fig. 12b).

The corresponding optimal phase-1 and phase-2 log-D-criteria $\Phi$ and CPU times are gathered in Table 2. The CPU times for Sobol(6) and Factorial(6) correspond to the total execution time of all 5 multistart runs.

The table shows that WDA(Grid(90)) and wMaxVol (Grid(90)) produce the global minimizer $\xi^*$, while Sobol(6) and Factorial(6) remain stuck at the local minimizer $\hat{\xi}$. The difference of computation times between wMaxVol and WDA is similar to the single start case. As expected, the multistart approach for Factorial(6) produces a lower objective value than in the single start case, but at a significantly higher computational cost. Overall, the cost of multistart is substantial in comparison with the discretization strategies.

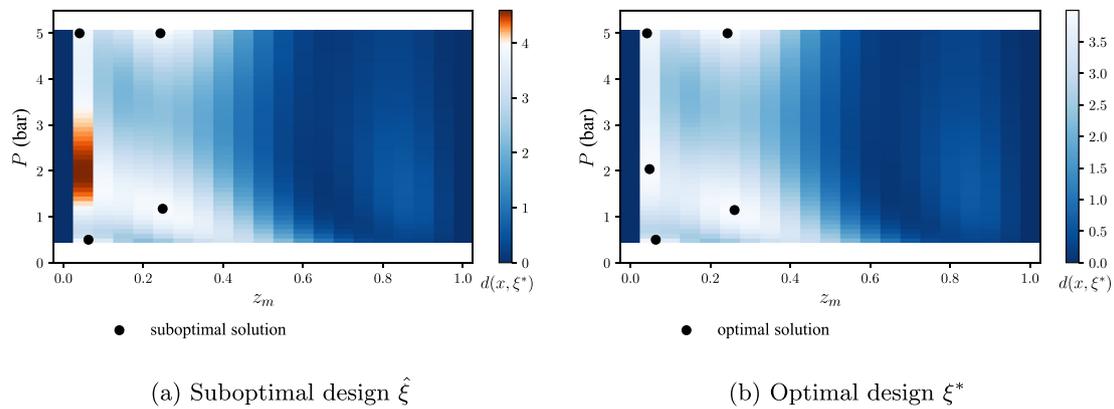

**Fig. 13.** Verification of phase-2 designs for the flash distillation (grid of 966 points).

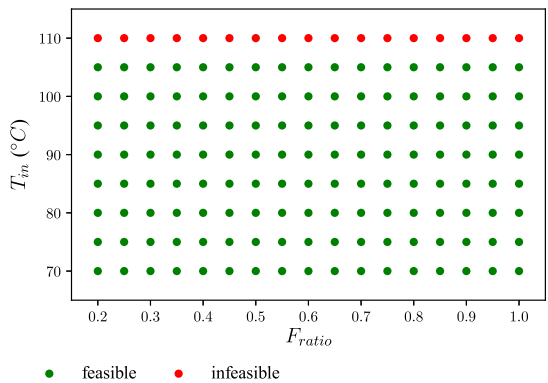

**Fig. 14.** Status of simulations for the tubular reactor.

**Table 2**
Optimal log-D-criteria and CPU times for the flash distillation ($N = 6$) with multi-start.

| Strategy | | Φ | CPU time (s) |
|---|---|---|---|
| `Sobol(6)` | Phase 1 | −4.839 | 97.29 |
| | Phase 2 | −7.928 | |
| `Factorial(6)` | Phase 1 | −3.373 | 112.34 |
| | Phase 2 | −7.928 | |
| `WDA(Grid(90))` | Phase 1 | −7.558 | 1.19 + 153 |
| | Phase 2 | **−7.935** | 16.79 |
| `wMaxVol(Grid(90))` | Phase 1 | −7.6097 | 1.19 + 0.421 |
| | Phase 2 | **−7.935** | 16.79 |

In order to assess the optimality of the designs $\hat{\xi}$ (suboptimal) and $\xi^*$ (putative globally optimal) on the continuous domain, we perform a verification on a (close to continuous) finer grid ($21 \times 46 = 966$ points). Fig. 13a exhibits $d$ values above $P = 4$ in the vicinity of $(z_m, P) = (0.048, 2.039)$, an experiment that exists in $\xi^*$ but not in $\hat{\xi}$. This proves that $\hat{\xi}$ is not globally optimal. Fig. 13b suggests that $\xi^*$ is likely to be the global minimizer, as the maximum of $x \mapsto d(x, \xi^*)$ on the grid is achieved at the support points of $\xi^*$ with a value of $P = 4$.

### 6.3. Tubular reactor

The tubular reactor is an arduous test problem ; some configurations (when $T_{in} > 105$) do not satisfy the constraints of the design at equilibrium. Fig. 14 illustrates the status (feasible or infeasible) of simulations generated on a $17 \times 9$ grid. Note that the feasible set depends on the current estimate of the model parameters $p$.

Although NLP solvers may converge towards a feasible point starting from an infeasible initial point, we require the set of initial phase-1 points to be feasible. Initial points for `Sobol` and `Factorial`, and candidate experiments for `WDA` and `wMaxVol` are therefore generated within the feasible domain $[0.2, 1] \times [70, 105]$. The `Factorial` strategy generates experiments in subdomains of the domain in a recursive fashion, such that the volume of each subdomain decreases by 20% at each iteration.

#### 6.3.1. Single start

Fig. 15a and b illustrate the phase-1 and phase-2 solutions for `Sobol` and `Factorial` with $N = 2$ experiments. Their phase-

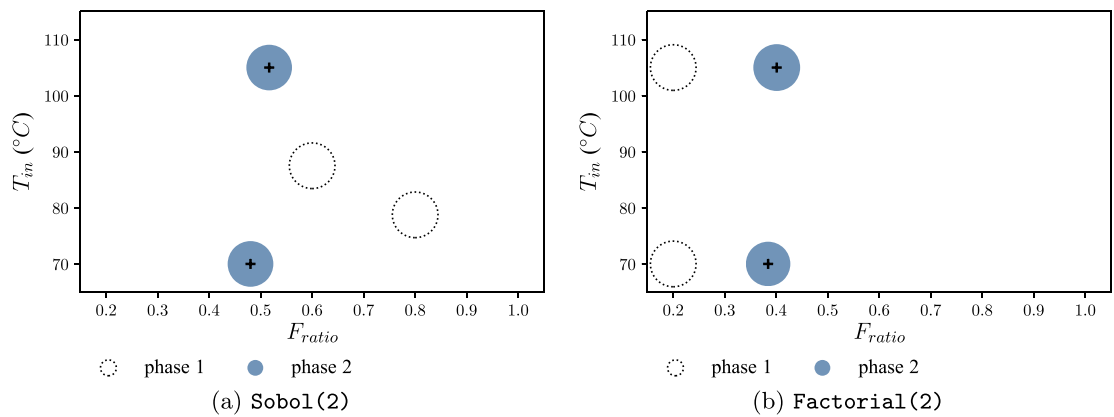

**Fig. 15.** `Sobol` and `Factorial` designs for the tubular reactor with 2 experiments.

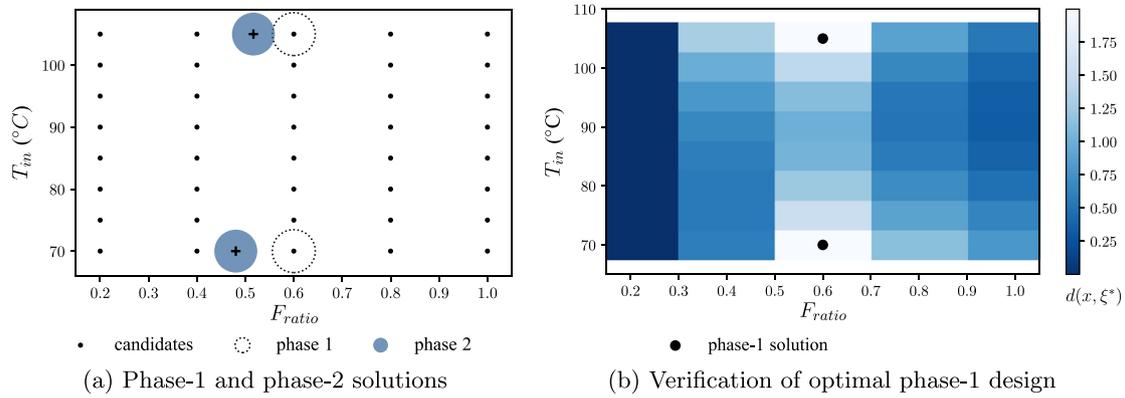

(a) Phase-1 and phase-2 solutions     (b) Verification of optimal phase-1 design

**Fig. 16.** WDA and wMaxVol designs on Grid(40) for the tubular reactor with 40 candidates.

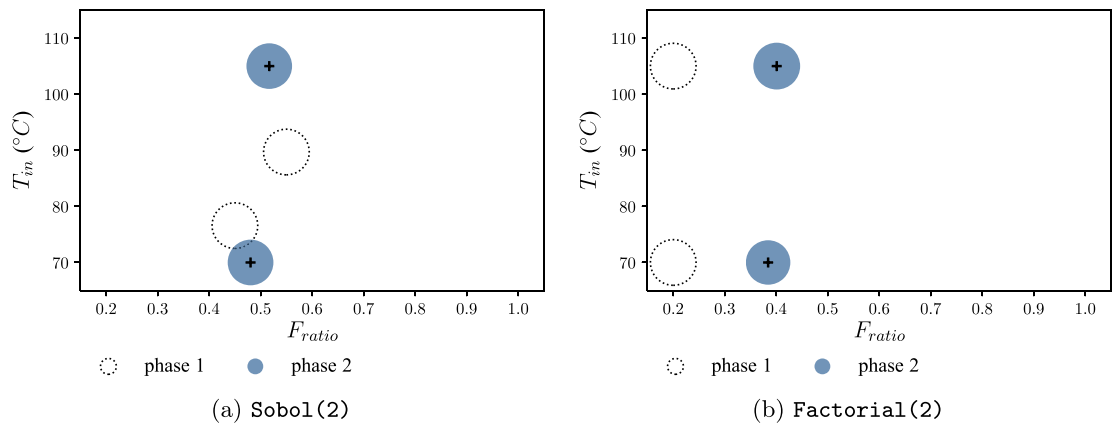

(a) Sobol(2)     (b) Factorial(2)

**Fig. 17.** Sobol and Factorial designs for the tubular reactor ($N = 2$) with 5 multistart runs.

**Table 3**
Optimal log-D-criteria and CPU times for the tubular reactor ($N = 2$) with single start.

| Strategy | | $\Phi$ | CPU time (s) |
| --- | --- | --- | --- |
| Sobol(2) | Phase 1 | 1.326 | 1.6 |
| | Phase 2 | **−0.633** | |
| Factorial(2) | Phase 1 | 6.710 | 1.51 |
| | Phase 2 | 0.275 | |
| WDA(Grid(40)) | Phase 1 | −0.324 | 0.97 + 0.55 |
| | Phase 2 | **−0.633** | 0.78 |
| wMaxVol(Grid(40)) | Phase 1 | −0.324 | 0.97 + 0.32 |
| | Phase 2 | **−0.633** | 0.78 |

**Table 4**
Optimal log-D-criteria and CPU times for the tubular reactor ($N = 2$) with multistart.

| Strategy | | $\Phi$ | CPU time (s) |
| --- | --- | --- | --- |
| Sobol(2) | Phase 1 | 0.050 | 7.81 |
| | Phase 2 | **−0.633** | |
| Factorial(2) | Phase 1 | 6.710 | 5.01 |
| | Phase 2 | −0.304 | |
| WDA(Grid(136)) | Phase 1 | −0.602 | 1.59 + 0.95 |
| | Phase 2 | **−0.633** | 0.49 |
| wMaxVol(Grid(136)) | Phase 1 | −0.602 | 1.59 + 0.75 |
| | Phase 2 | **−0.633** | 0.49 |

2 solutions are $\xi^* = \begin{Bmatrix} 0.517 & 0.480 \\ 105 & 70 \\ 0.499 & 0.501 \end{Bmatrix}$ and $\begin{Bmatrix} 0.384 & 0.401 \\ 70 & 105 \\ 0.474 & 0.526 \end{Bmatrix}$, respectively.

The identical results of WDA and wMaxVol on a coarse grid of 40 feasible candidate experiments are shown in Fig. 16a. Both strategies produce the phase-2 design $\xi^*$. The verification of optimality is shown in Fig. 16b: the phase-1 solution is optimal for the discrete domain Grid(40).

The corresponding optimal phase-1 and phase-2 log-D-criteria $\Phi$ and CPU times are gathered in Table 3. For the discretization strategies, the CPU times include i) the generation of Jacobians at the candidate experiments (0.97s) and ii) the convergence time of phase 1. Factorial(2) finds a local minimizer, while the other strategies produce the putative global minimizer $\xi^*$. Again, wMaxVol performs faster than WDA, albeit with a marginal improvement.

*6.3.2. Multistart*

Fig. 17a and b illustrate the phase-1 and phase-2 solutions of Sobol and Factorial for the best run among 5 multistart runs. Similarly to the single start case, Factorial(2) converges towards a local minimizer.

Fig. 18a illustrates the phase-1 and phase-2 solutions of the discretization strategies on a fine grid of 136 candidate experiments. WDA and wMaxVol produce the same phase-1 solution with 2 nonzero weights, very close to the putative global minimizer $\xi^*$, which is subsequently reached during phase 2. This suggests that the granularity of the grid with 40 candidate experiments is sufficient to find the global optimum. The global optimality of the phase-1 solution is certified on the discrete domain (Fig. 18b).

The corresponding optimal phase-1 and phase-2 log-D-criteria $\Phi$ and CPU times are gathered in Table 4. The CPU times for Sobol(2) and Factorial(2) correspond to the total execution time of all 5 multistart runs. For the discretization strategies, the CPU times include i) the generation of Jacobians at the candidate experiments (1.59s) and ii) the convergence time of phase 1. Using multistart, the putative global minimizer $\xi^*$ is reached with the Sobol(2), WDA(Grid(136)) and wMaxVol(Grid(136)) strategies. wMaxVol is again slightly faster than WDA.

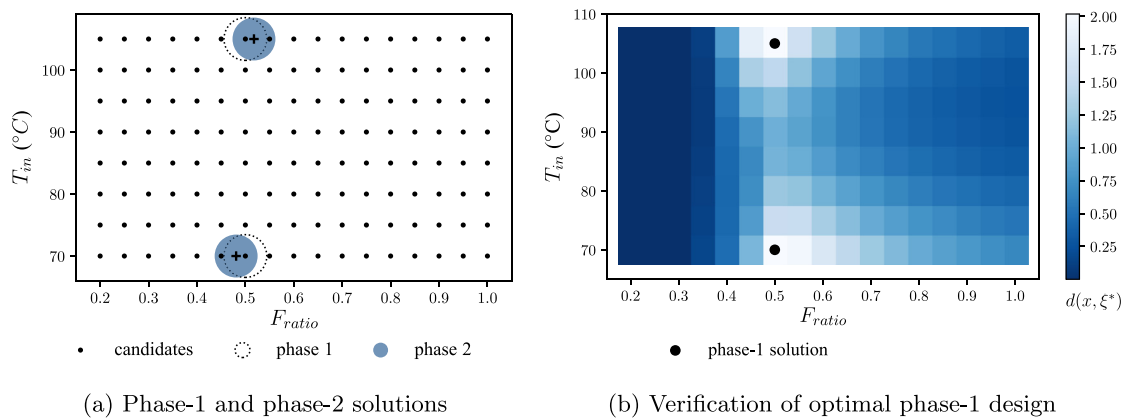

(a) Phase-1 and phase-2 solutions

(b) Verification of optimal phase-1 design

**Fig. 18.** WDA and wMaxVol designs on Grid(136) for the tubular reactor with 136 candidates.

## 7. Conclusions

The experimental design problem – the selection of optimal experiments and corresponding weights according to a given statistical criterion – is an arduous optimization problem on account of its nonconvexity. Gradient-based optimization solvers cannot guarantee the global optimality of the solution. Furthermore, a poor initial point may dramatically hinder the convergence of the method or the quality of the local optimum.

In this paper, we described a two-phase strategy to initialize and solve the nonconvex experimental design problem. We introduced WDA and wMaxVol, two problem-dependent discretization strategies that determine the weights of the relevant experiments over a fixed set of candidates, as well as an approximation of the optimal number of experiments. We assessed their performance against two pattern-based problem-independent strategies, Sobol and Factorial. On the one hand, the two-phase approach using discretization strategies only requires two optimization runs: a small convex phase-1 problem and a nonconvex phase-2 problem initialized with a satisfactory initial point. This reduces the effort required to solve the experimental design problem to global optimality. On the other hand, Sobol and Factorial are generic and easy to implement, and may produce good results in combination with multistart (independent runs starting from different initial points), albeit at a higher cost.

When the set of candidate experiments is made finer, the discretization strategies perform consistently better than the pattern-based strategies and produce phase-1 solutions that are close to the optimal phase-2 solutions, which avoids numerous solver iterations. It is therefore worthwhile to successively refine the discretization in the neighborhood of the solution, in order to iteratively improve the initial guess. The pattern-based strategies suffer from one major drawback: the number of relevant experiments is not known a priori, while it is automatically inferred by the discretization strategies for a given discretization. WDA and wMaxVol thus answer the key interrogations in experimental design: how many experiments should be picked, and which ones.

Parametric uncertainty is of decisive importance, since it may impact the feasibility of the physical constraints, the number of optimal experiments and the optimal design. In this paper, the model parameters were fixed to their current estimates, and parametric uncertainty was handled by our iterative model validation and adjustment framework. In future research, we plan to extend the methods described in this paper to solve models with explicit uncertain parameters.

## Declaration of Competing Interest

The authors declare that they have no known competing financial interests or personal relationships that could have appeared to influence the work reported in this paper.

## Acknowledgements

The work of Charlie Vanaret, Philipp Seufert, Jan Schwientek and Michael Bortz was developed in the Fraunhofer Cluster of Excellence "Cognitive Internet Technologies". The work of Gleb Karpov, Gleb Ryzhakov and Ivan Oseledets was supported by the Ministry of Science and Higher Education of the Russian Federation Agreement 075-10-2020-091 (grant 14.756.31.0001).

## Appendix A. wMaxVol on univariate problems: a sequential design generation algorithm

According to Fedorov and Leonov (2014), Wynn (1970), one of the methods to construct a D-optimal design is to sequentially add points to the design $\xi^{(k)} = \left\{ \begin{array}{c} supp(\xi^{(k)}) \\ w^{(k)} \end{array} \right\}$ such that:

$$x^{(k+1)} = \arg\max_{x \in \mathcal{D}} d(x, \xi^{(k)}) \qquad (A.1)$$

The new design $\xi^{(k+1)}$ can be described as follows:

- if $x^{(k+1)} \in supp(\xi^{(k)})$, the set of support points remains unchanged, all the weights are multiplied by $(1 - \alpha^{(k)})$ and $\alpha^{(k)}$ is added to the weight of $x^{(k+1)}$ ;
- if $x^{(k+1)} \notin supp(\xi^{(k)})$, $x^{(k+1)}$ is added with weight $\alpha^{(k)}$ to the set of support points and all the other weights are multiplied by $(1 - \alpha^{(k)})$ ;

where the step $\alpha^{(k)} \in (0, 1)$ may be picked according to various strategies, for instance such that:

$$\sum_{k=0}^{+\infty} \alpha^{(k)} = +\infty \quad \text{and} \quad \lim_{k \to +\infty} \alpha^{(k)} = 0. \qquad (A.2)$$

In the following, we consider the case of univariate observations ($m = 1$) and write the blocks $A_i$ and $C_i$ as row vectors $a_i$ and $c_i$, respectively. We prove that wMaxVol satisfies the two conditions of the family of sequential design generation algorithms:

1. the next experiment $x_i$ maximizes the sensitivity function $x \mapsto d(x, \xi^{(k)})$ ;

2. there exists $\alpha^{(k)}$ compatible with Eq. (A.2).

The theory developed for the sequential design generation algorithms subsequently applies. At each iteration of Algorithm 2, an intermediate full-column rank matrix $\tilde{A}^{(k)}$ uniquely determines the corresponding unnormalized information matrix $\underline{\mathcal{I}}(\xi^{(k)})$. We adopt the notations from Fedorov and Leonov (2014) and introduce:

$$(\tilde{A}^{(k)})^T \tilde{A}^{(k)} = \underline{\mathcal{I}}(\xi^{(k)}) = s^{(k)} \mathcal{I}(\xi^{(k)}), \quad (A.3)$$

where $s^{(k)} = s^{(0)} + k$ is the current number of non-unique points in the design $\xi^{(k)}$ and $s^{(0)}$ is the number of points in the initial design of $\tilde{A}$.

The wMaxVol expansion step yields the augmented matrix $\tilde{A}^{(k+1)} = \begin{bmatrix} \tilde{A}^{(k)} \\ a_i \end{bmatrix}$. The corresponding information matrix $\underline{\mathcal{I}}(\xi^{(k+1)})$ can be written as:

$$\underline{\mathcal{I}}(\xi^{(k+1)}) = \underline{\mathcal{I}}(\xi^{(k)}) + a_i^T a_i. \quad (A.4)$$

Consider the $i$th row $c_i = a_i(\tilde{A}^{(k)})^+$ of $C = A(\tilde{A}^{(k)})^+$. Using the properties of the Moore-Penrose inverse and the assumption that $\tilde{A}^{(k)}$ is full rank, we get:

$$c_i c_i^T = a_i (\tilde{A}^{(k)})^+ (\tilde{A}^{(k)})^{+T} a_i^T = a_i \underline{\mathcal{I}}^{-1}(\xi^{(k)}) a_i^T. \quad (A.5)$$

The sensitivity function (Eq. (10)) is in this case written:

$$d(x_i, \xi^{(k)}) = a_i \mathcal{I}^{-1}(\xi^{(k)}) a_i^T. \quad (A.6)$$

where $x_i \in \mathcal{D}$ is the experiment corresponding to the row $a_i$ in the design matrix. Normalizing Eq. (A.5), we conclude:

$$c_i c_i^T = \frac{1}{s^{(k)}} d(x_i, \xi^{(k)}). \quad (A.7)$$

We thus established that picking the experiment $x_i$ with the maximal $\|c_i\|_2^2$ (Line 7 in Algorithm 2) amounts to picking the point with the maximal $d(x_i, \xi^{(k)})$.

Now let us derive the step size $\tilde{\alpha}^{(k)}$ for the wMaxVol expansion step (adding one point to the current set of $s^{(k)}$ points). $\mathcal{I}(\xi^{(k+1)})$ can be written as a convex combination:

$$\mathcal{I}(\xi^{(k+1)}) = (1 - \tilde{\alpha}^{(k)}) \mathcal{I}(\xi^{(k)}) + \tilde{\alpha}^{(k)} \mu(x_i). \quad (A.8)$$

Normalizing Eq. (A.4), we obtain:

$$\mathcal{I}(\xi^{(k+1)}) = \frac{s^{(k)}}{s^{(k)}+1} \mathcal{I}(\xi^{(k)}) + \frac{\mu(x_i)}{s^{(k)}+1} \quad (A.9)$$

From Eqs. (A.8) and (A.9), we finally conclude:

$$\tilde{\alpha}^{(k)} = \frac{1}{s^{(k)}+1} = \frac{1}{s^{(0)}+k+1}. \quad (A.10)$$

The harmonic series $\sum_{k=0}^{+\infty} \tilde{\alpha}^{(k)}$ diverges and $\lim_{k \to +\infty} \tilde{\alpha}^{(k)} = 0$, which is compatible with Eq. (A.2). Thus, we showed that wMaxVol satisfies the two conditions of the family of sequential design generation algorithms. From Theorem 3.2 in Fedorov and Leonov (2014), we state that if the initial design $\tilde{A}$ is regular, i.e. $\det(\tilde{A}^T \tilde{A}) \neq 0$, we have:

$$\lim_{k \to \infty} \Psi(\mathcal{I}(\xi^{(k)})) = \min_{\xi} \Psi(\mathcal{I}(\xi)). \quad (A.11)$$

### Appendix B. wMaxVol on multivariate problems: a statistical analysis

Our goal is to perform a statistical analysis of the convergence of wMaxVol in the multivariate case and prove that it is indeed suited for D-optimal experimental design. Sequential design generation algorithms iteratively seek an optimal design by either inserting a new point into an existing design or by correcting the weight of an existing point in the design. A precise estimation of the weights is of crucial importance: a design $\xi$ that contains the design points $supp(\xi^*)$ of the optimal design $\xi^*$, however paired with arbitrary weights, is detected as suboptimal by the Kiefer-Wolfowitz equivalence theorem, that we recall here:

$$\begin{cases} d(x, \xi^*) \leq P, & \forall x \in \mathcal{D} \\ d(x_i^*, \xi^*) = P, & \forall x_i^* \in supp(\xi^*) \end{cases} \quad (B.1)$$

Sequential algorithms usually spend numerous iterations seeking the optimal weights associated with the optimal design points, therefore they produce suboptimal designs until the values of the weights settle down. wMaxVol computes the weight of a design point as a rational approximation (line 13 in Algorithm 2) of the optimal weight. In addition, numerical errors may affect the sensitivity value of the final design: it may not be strictly equal to the optimal bound $P$, as required by Eq. (B.1). We therefore introduce two metrics that rely on Eq. (B.1) to measure the progress of wMaxVol and its convergence towards a global minimum. If wMaxVol indeed converges towards a D-optimal design, we expect that the two metrics tend to 0 when $k$ goes to $+\infty$.

The first metric is based on the equality $d(x_i^*, \xi^*) = P$, $\forall x_i^* \in supp(\xi^*)$ of Eq. (B.1). It measures the variance of the sensitivity function around $P$ and penalizes its deviation from $P$:

$$\Delta(\xi^{(k)}) = \max(\overline{d}(\xi^{(k)}), P) - \min(\underline{d}(\xi^{(k)}), P). \quad (B.2)$$

where

$$\overline{d}(\xi^{(k)}) = \max_x d(x, \xi^{(k)}), \quad x \in supp(\xi^{(k)})$$

$$\underline{d}(\xi^{(k)}) = \min_x d(x, \xi^{(k)}), \quad x \in supp(\xi^{(k)}).$$

The second metric is based on the inequality $d(x, \xi^*) \leq P$, $\forall x \in \mathcal{D}$ of Eq. (B.1). It measures the fraction of points $x$ outside the current design $\xi^{(k)}$ whose sensitivity values exceed $P$:

$$q(\xi^{(k)}) = \frac{|\{x \in \mathcal{D} \setminus supp(\xi^{(k)}) \mid d(x, \xi^{(k)}) > P\}|}{|\mathcal{D} \setminus supp(\xi^{(k)})|}. \quad (B.3)$$

In early iterations, when the current design $\xi^{(k)}$ is likely to be suboptimal, the sensitivity values of points outside $\xi^{(k)}$ may exceed $P$. However, since the values of the weights in $\xi^{(k)}$ are iteratively refined, the number of points outside $\xi^{(k)}$ whose sensitivity values exceed $P$ must tend to zero for the equivalence theorem to be verified.

We now introduce an example on which we derive a statistical study. Consider the experimental design problem with a linear model and a multivariate output. Let $\{\phi_1, \ldots, \phi_P\}$ be a set of basis functions and $\phi_i$, $i \in \{1, \ldots, P\}$, can be factorized into univariate polynomials $T_{\lambda_{ij}}$ of degree $\lambda_{ij}$:

$$\phi_i(x) = \prod_{j=1}^{n} T_{\lambda_{ij}}(x_j),$$

$\lambda_i \in \mathbb{N}^n$ are constructed according to a special case of hyperbolic truncation of polynomial expansion: $\|\lambda_i\|_1 \leq d$, where $d$ is the maximum degree of the basis functions' polynomials (see Blatman and Sudret (2011)). We chose Chebyshev polynomials defined by:

$$\begin{cases} T_0(x) = 1 \\ T_1(x) = x \\ T_{i+1}(x) = 2xT_i(x) - T_{i-1}(x), \quad \forall i \geq 1 \end{cases}$$

The $i$th column of the matrix $\varphi(x)^T \in \mathbb{R}^{m \times P}$ is formed by the value of the basis function $\phi_i$ and its partial derivatives:

$$\varphi(x)^T = \begin{pmatrix} \phi_1(x) & \cdots & \phi_P(x) \\ \frac{\partial \phi_1(x)}{\partial x_1} & \cdots & \frac{\partial \phi_P(x)}{\partial x_1} \\ \vdots & \ddots & \vdots \\ \frac{\partial \phi_1(x)}{\partial x_n} & \cdots & \frac{\partial \phi_P(x)}{\partial x_n} \end{pmatrix}. \quad (B.4)$$

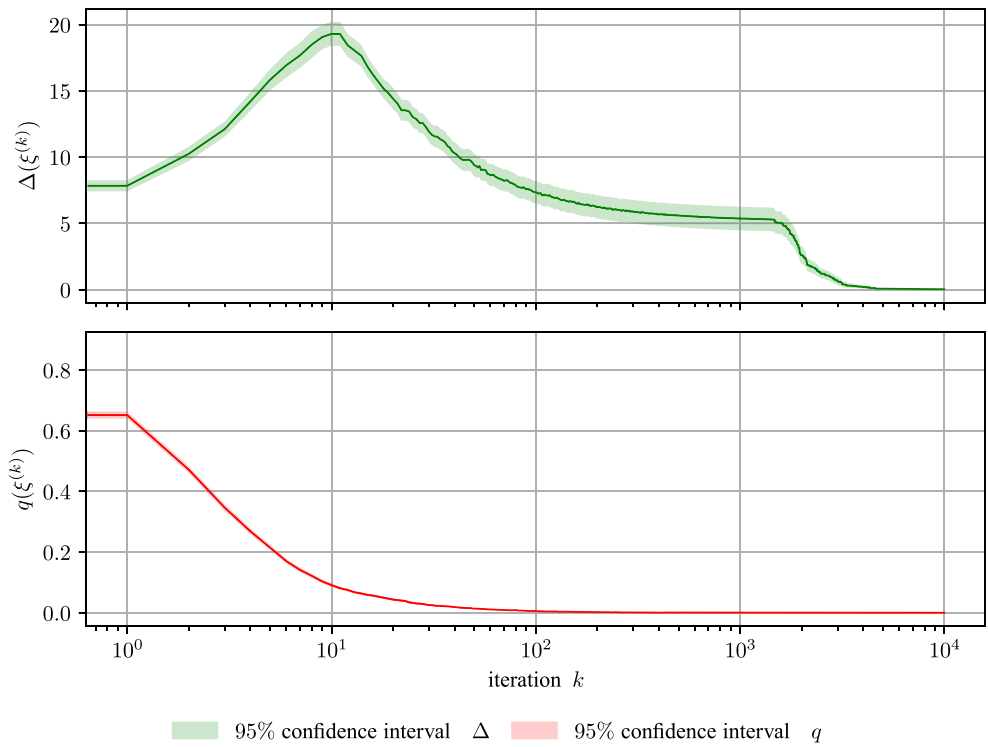

Fig. B.19. Convergence of `wMaxVol` on Chebyshev polynomials for $n = 5$.

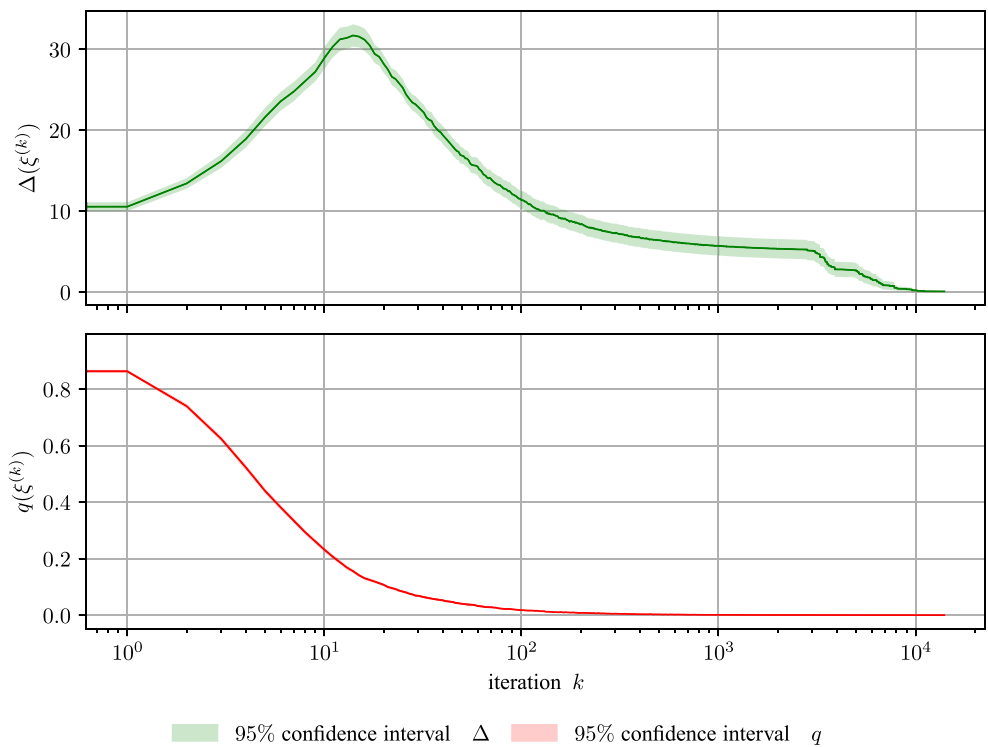

Fig. B.20. Convergence of `wMaxVol` on Chebyshev polynomials for $n = 7$.

We ran `wMaxVol` on two instances of the problem for $n = 5$ and $n = 7$. The values of the parameters are given in Table B.5. For each of the 100 runs, each of the 200 candidate experiments is randomly generated in $[-1, 1]^n$ using a Latin hypercube strategy.

The two metrics $\Delta$ and $q$ are reported in Fig. B.19 (for $n = 5$) and B.20 (for $n = 7$). The $x$ axis represents the current iteration number $k$ in logarithmic scale. For both metrics, the mean (solid line) over the 100 runs and a 95% confidence interval (envelope) are shown. As expected, both metrics start with a positive value and tend to 0 when the number of iterations increases. This demonstrates that the global minimum is attained. This statistical study suggests that `wMaxVol` is indeed suited for multivariate experimental design.

**Table B.5**
Parameters of the `wMaxVol` statistical study.

| Parameter | Instance 1 | Instance 2 |
|---|---|---|
| $n$ design variables | 5 | 7 |
| $K$ candidate experiments | 200 | 200 |
| $P$ model parameters | 36 | 64 |
| runs | 100 | 100 |
| iterations | 10,000 | 14,000 |